\documentclass{article}

\usepackage{amsfonts,amsmath, amsthm, amscd, amssymb, graphicx, color, mathrsfs, stmaryrd}

\def\findemo{\hfill \rule{6pt}{6pt}}

\begin{document}

\newtheorem{thm}{Theorem}[section]
\newtheorem{theo}[thm]{Theorem}
\newtheorem{prop}[thm]{Proposition}
\newtheorem{coro}[thm]{Corollary}
\newtheorem{lema}[thm]{Lemma}
\newtheorem{defi}[thm]{Definition}
\newtheorem{ejem}[thm]{Example}
\newtheorem{rema}[thm]{Remark}
\newtheorem{fact}[thm]{Fact}
\newtheorem{open}[thm]{PROBLEM}

\newcommand{\lip}{{\rm Lip}}
\newcommand{\lipo}{ {\rm Lip}_0 }
\newcommand{\conv}{\operatorname{conv}}
\newcommand{\cconv}{\overline{\conv}}

\title{Covering $B_X$ by finitely many convex sets}
\author{M. Raja\thanks{This research has been supported by:  Fundaci\'on S\'eneca -- ACyT Regi\'on de Murcia, project 21955/PI/22; and grant PID2021-122126NB-C32  funded by MCIN/AEI/ 10.13039/501100011033 and by “ERDF A way of making Europe”.}}
\date{February, 2025}

\maketitle

\begin{abstract}
Given a finite covering by closed convex sets of $B_X$, the unit ball of an infinite-dimensional Banach space, we investigate whether there is a set of the covering that contains balls of radius close to $1$ and (a) arbitrarily high finite dimension or (b) infinite dimension. In case (a) the answer is affirmative, but for the case (b) we just get radius close to $1/2$ and finite codimension under much more restrictive hypotheses.
\end{abstract}




\section{Motivation}

Baire's theorem implies that in any countable cover of the unit ball $B_X$ of a Banach space $X$ by closed sets there is, at least, one set from the cover that contains a smaller ball of positive radius. It seems quite natural to expect a little more if the cover by closed sets is finite, or the sets are moreover convex. However, very simple examples show that we will not have balls of bigger radius being finite the cover of $B_X$.\\

On the other hand, a well known result of Lyusternik and Shnirel’man, see \cite{LS, Matusek}, that any cover of ${\Bbb S}^{n-1}$ (the Euclidean sphere in ${\Bbb R}^n$) by not more than $n$ closed sets, there is one of them containing a pair of antipodal points. The Lyusternik-Shnirel’man theorem easily implies the following result: if $X$ is an infinite dimensional Banach space, then in any finite cover of $B_X$ by closed convex sets, there is one of them that contains a diameter, in other words, a $1$-dimensional ball of radius $1$.\\

Our aim in this paper is to maximize the radius and the (vectorial) dimension of the balls contained in the sets of a finite cover of $B_X$ by closed convex sets. Let us state our basic definition. We say that a set $A \subset X$ contains an $n$-dimensional (respectively, infinite dimensional, finite codimensional) ball of radius $\lambda>0$ if there is $x \in A$ and a subspace $Y \subset$ of dimension $n$ (resp. infinite dimension, finite codimension) such that
$$ x + \lambda B_Y \subset A, $$
where $B_Y=B_X \cap Y$. We will see that it is possible to find balls of arbitrarily high dimension and radius close to $1$ in finite covers of $B_X$ by closed convex sets, but infinitely dimensional balls seems to be scarcer.\\

The organization of the paper is the following. In the next section we exhibit two relevant examples and provide some auxiliary lemmas. Then, in section 3, we address the problem of finding high-dimensional balls of large radius. It turns out that the solution follows easily from some not well known early results of Milman. However, for the sake of completeness we provide a full proof based in  Milman's concentration of measure phenomenon and Dvoretzky's theorem. The aim in the last section is to find infinite-dimensional (actually, of finite codimension) balls of large radius in a set that covers $B_X$ with finitely many of its translates. The idea here is to find a projection onto a suitable finite codimensional subspace that approximatively preserves the set.\\

We consider that the results in this paper are far from being optimal, but they point to a promising line of research.
All the Banach spaces considered are real. Our notation is totally standard and we address to generic references for any unexplained definition \cite{JL, banach}.

\section{Convex decompositions of $B_X$}

When we started this work, the aim was to find finite codimensional balls of large enough radius in the pieces of a finite decomposition of the unit ball $B_X$ of an infinite dimensional Banach space $X$, mainly motivated by its relation to {\it asymptotic uniformly smooth renormings}, 
see \cite{raja}. We will see that the answer depends largely on the geometry of $X$.\\

Firstly, we will show with a simple example that the largest finite codimensional balls contained in the pieces of a convex cover of the unit ball of the Hilbert space $H$ may have small radius. 

\begin{ejem}
For every $k \in {\Bbb N}$ it is possible to decompose $B_H$ as the union of $2k$ congruent closed convex sets such that if $r_k$ is the supremum of the radius of finite-codimensional balls contained in any of those sets, then $r_1<1$ and $\lim_k r_k=0$.
\end{ejem}

\noindent
{\bf Proof.}
We will take $H=\ell_2$ with the coordinates indexed by ${\Bbb N} \cup \{0\}$. Fix $k \geq 1$ and
for $1 \leq j \leq 2k$ consider the set
$$ A_j = \{ (x_n) \in B_{\ell_2}: (-1)^j x_0 \leq 1/2 - k \sum_{m=0}^\infty x_{2km+j}^2 \} .$$
These sets are closed and convex, and moreover, $\bigcup_{j=1}^{2k} A_j=B_{\ell_2}$. 
Indeed, assume that some $x =(x_n) \in B_{\ell_2}$  belongs to none of the sets. Then we have the reversed inequalities 
$$ (-1)^j x_0 > 1/2 - k \sum_{m=0}^\infty x_{2km+j}^2 ,$$
whose sum over $1 \leq j \leq 2k$ gives
$$ 0 > k -k \sum_{n=1}^\infty x_{n}^2 .$$
We deduce
$$  1 <   \, \sum_{n=1}^\infty x_{n}^2 \,  \leq \|x\|^2  $$
that is a contradiction with $x \in B_{\ell_2}$.
Now we will estimate the radius of a finite codimensional ball contained in $A_j$. 
Let $Y$ be a finite codimensional subspace. Obviously, $Y$ meets the vectors supported on coordinates of the form $2km+j$ in a infinite-dimensional subspace of $\ell_2$, say $Z$. 
We can choose a $2$-dimensional section given by the direction of $x_0$ and vector from $Z$, and such that the coordinates of the center of the ball, except $x_0$, are negligible. This section reduces the problem to dimension $2$, where some Euclidean gymnastics easily provides this upper bound for $r_k$
$$ r_k \leq \sqrt{  \frac{1}{2k}\left(1 - \frac{1}{k} \right)  + \frac{1}{2k} \sqrt{ \left(1 - \frac{1}{k} \right)^2 + 3} }. $$
The claimed properties of $(r_k)$ follow easily.\findemo\\

We do not know whether it is possible to produce a similar example in $H$ with infinite dimensional balls instead of finite codimensional.\\

Now we will show that $c_0$ behaves very differently from Hilbert space when it comes to finite convex decompositions of the unit ball.

\begin{ejem}\label{cecero}
Let $A_j \subset B_{c_0}$ for $j=1,\dots,k$ be convex closed subsets such that
$B_{c_0} = \bigcup_{j=1}^k A_j$. Then for some $j$, there is a finite codimensional ball of radius $1$ contained in $A_j$.
\end{ejem}

\noindent
{\bf Proof.} 
Indeed, consider $\overline{A_j}^{w^*}$  in $\ell_\infty$. Note that the closure only adds points which are in 
$\ell_\infty \setminus c_0$ since $A_i$ is weakly closed. The extreme points of  $B_{\ell_\infty}$ is $\{-1,1\}^{\Bbb N}$, that becomes the Cantor set when considered with the weak$^*$ topology.
Then $B_{\ell_\infty} \subset \bigcup_{j=1}^k \overline{A_j}^{w^*}$ and, in particular,  
$\{-1,1\}^{\Bbb N}$ is covered. The Baire category theorem (one could argue whether for finitely many sets Baire name is adequate) applied to $\{-1,1\}^{\Bbb N}$ provides some $\overline{A_j}^{w^*}$ that contains a $w^*$-clopen set of the form
$$ (a_1,\dots,a_i) \times \{-1,1\}^{\Bbb N} .$$
The $w^*$-closed hull of this set is contained by  $\overline{A_j}^{w^*}$ and therefore 
$$(a_1,\dots,a_i) \times [-1,1]^{\Bbb N} \subset \overline{A_j}^{w^*}$$
 with $x_0 = (a_1,\dots,a_i,0,0,\dots)$. As $x_0 \in c_0$ we have 
$x_0 + B_Y \subset A_k$ being $Y$ the subspace of $c_0$ having null the first $i$ coordinates.\findemo\\

We will finish with some observations that will be useful in next sections.

\begin{lema}\label{tech_1}
Let $A \subset X$ be bounded convex with nonempty interior and assume $0 \in A^{\mathrm{o}}$. 
Then, for every $\varepsilon>0$ there is $\delta>0$ such that
$$ A + \delta B_X \subset (1+\varepsilon) A, ~~\mbox{and}~~ (1+\delta) A  \subset A + \varepsilon B_X. $$
\end{lema}

\begin{lema}\label{tech_2}
Let $A \subset X$ be a closed set with nonempty interior. If $A$ is covered by finitely many convex closed sets 
$(A_j)_{j=1}^k$, then $A$ is covered by those among 
$(A_j)_{j=1}^k$ with nonempty interior. 
\end{lema}

\noindent
{\bf Proof.}
By Baire's theorem, the set $ \bigcup_{j=1}^k A_j^{\mathrm{o}}$ is dense in $A$, thus
$$ A \subset \overline{ \bigcup_{j=1}^k A_j^{\mathrm{o}}} =   \bigcup_{j=1}^k  \overline{ A_j^{\mathrm{o}}} = 
\bigcup \{ A_j:  A_j^{\mathrm{o}} \not = \emptyset \}, $$
as claimed.\findemo

\section{Arbitrarily high-dimensional Euclidean balls}

Let us start by reminding the following application of the concentration of measure phenomenon for functions on the Euclidean sphere proved by Milman in \cite{Milman1}. We will follow the formulation given in \cite[Proposition 12.3]{BL}.

\begin{lema}\label{lema_BL}
Let $f$ be a real-valued function on ${\Bbb S}^{N-1}$ which is $\tau$-Lipschitz with respect to the Euclidean norm, and let 
$0<\varepsilon<\tau/2$. Then there is an $n$-dimensional subspace $F$ of ${\Bbb R}^N$ with $n \geq CN\varepsilon^2/|\log (\varepsilon/\tau)|\tau^2$ (where $C$ is universal constant), and there is a number $\lambda_0$ such that
$$ \lambda_0 - \varepsilon \leq f(x) \leq   \lambda_0 + \varepsilon  $$
for every $x \in F \cap {\Bbb S}^{N-1}$.
\end{lema}

The set of values where a continuous function $f$ on the sphere are approximative concentrated on arbitrarily high finite dimensional spaces is called the {\it spectrum} of $f$ by Milman. Lemma \ref{lema_BL} establishes that the spectrum is nonempty, see \cite{Milman3} for an account of the development of this theory.\\

The next three results are also due to Milman, however their traceability was not easy to me because, to my knowledge, they do not appear in books despite their great importance. For instance, the infinite dimensional version of Lemma \ref{lema_2} appears in \cite{gromov} and \cite{Pestov} with no proof. For that reason, we will include qualitative proofs using Lemma \ref{lema_BL} and Dvoretzky's theorem \cite[Theorem 12.10]{BL}.

\begin{lema}\label{lema_1}
Let $n \in {\Bbb N}$ and $\varepsilon \in (0,1/2)$. There is $N(n, \varepsilon)$ such that for all $N \geq N(n,\varepsilon)$ and $A_1, A_2$ a covering of ${\Bbb S}^{N-1}$ by two closed sets, there is an $n$-dimensional subspace $F$ of ${\Bbb R}^N$ and some $i \in \{1,2\}$. such that
$$ F \cap {\Bbb S}^{N-1} \subset A_i +\varepsilon B_{ {\Bbb R}^N} .$$
\end{lema}

\noindent
{\bf Proof.}
Take $f(x)=d(A_1, x)$, that is $1$-Lipschitz and apply Lemma~\ref{lema_BL} with $\varepsilon/2$ so there is $\lambda_0$ and $F \subset B_{ {\Bbb R}^N}$ an $n$-dimensional subspace such that
$$ \lambda_0 - \varepsilon/2 \leq f(x) \leq   \lambda_0 + \varepsilon/2  $$
for $x \in F \cap {\Bbb S}^{N-1}$ . Then either $\lambda_0 \leq \varepsilon/2$ or $\lambda_0 > \varepsilon/2$. In the first case we have 
$$ 0 \leq f(x) \leq \varepsilon $$
for every $x \in F \cap {\Bbb S}^{N-1}$
and therefore $A_1$ does the work. On the other hand, if $\lambda_0 > \varepsilon/2$, then $f(x) >0$ for for every $x \in F \cap {\Bbb S}^{N-1}$. Consequently,
$$ x \in A_1^c \subset A_2 \subset A_2 +  \varepsilon B_{ {\Bbb R}^N} $$
that finishes the proof.\findemo\\

Our Lemma \ref{lema_1} is \cite[Corollary 3]{Milman1} without the precise estimation of the dimension given by Milman.

\begin{lema}\label{lema_2}
Let $n, k \in {\Bbb N}$ and $\varepsilon \in (0,1/2)$. There exists $N(n, k, \varepsilon)$ such that for all $N \geq N(n, k,\varepsilon)$ and $(A_i)_{i=1}^k$ a covering of ${\Bbb S}^{N-1}$ by closed sets, there is an $n$-dimensional subspace $F$ of ${\Bbb R}^N$ and some $i \in \{1, \dots, k\}$ such that 
$$ F \cap {\Bbb S}^{N-1} \subset A_i +\varepsilon B_{ {\Bbb R}^N} . $$
\end{lema}

\noindent
{\bf Proof.}
The case $k=2$ is just the previous lemma. Assume that ${\Bbb S}^{N-1}$ is covered by three sets $A_1, A_2, A_3$ and apply Lemma \ref{lema_1} to the cover $A_1$, $A=A_2 \cup A_3$ with $\varepsilon/2$ and
$$ N \geq N(N(n,\varepsilon/2), \varepsilon/2). $$ 
Then we can find $E$ of dimension at least $N(n,\varepsilon/2)$ such that either 
$$A_1+(\varepsilon/2)B_{ {\Bbb R}^N} ~~\mbox{or}~~ A+(\varepsilon/2)B_{ {\Bbb R}^N}$$ 
contains $S_E$. In the first case, we are done. In the second one, we have a covering of $S_E$ by two closed sets 
$A_2+(\varepsilon/2)B_{ {\Bbb R}^N}$ and $A_3+(\varepsilon/2)B_{ {\Bbb R}^N}$, so we can apply again Lemma \ref{lema_1} in order to finish the $3$-set case. For a general $k$, take 
$$ N \geq N(n, k, \varepsilon) := N(N(\dots N(n,\varepsilon/k),\varepsilon/k),\varepsilon/k), $$ 
where the function $N(n,\varepsilon/k)$ is composed $k-1$ times with itself.\findemo

\begin{theo}\label{main_teo}
Let $n, k \in {\Bbb N}$ and $\varepsilon \in (0,1/2)$. There exists ${\mathcal N}(n, k, \varepsilon)$ such that for all the normed spaces $X$ of dimension $N \geq {\mathcal N}(n, k,\varepsilon)$ and $(A_i)_{i=1}^k$ a covering of $S_X$ by closed sets, there is an $n$-dimensional subspace $F$ of $X$ and some $i \in \{1, \dots, k\}$ such that 
$$ F \cap S_X \subset A_i +\varepsilon B_{X} . $$
\end{theo}

\noindent
{\bf Proof.} By Dvoretzky's theorem, see \cite[Theorem 12.10]{BL} for instance, if the dimension of $X$ is large enough, we can find a Euclidean sphere 
${\Bbb S}$ between $S_X$ and $(1+\varepsilon/3)S_X$ when restricted to some subspace $E$ of dimension $N(n, k, \varepsilon/3)$, the number given by Lemma \ref{lema_2}. If $(A_i)_{i=1}^k$ a covering of $S_X$, then 
$$ (A_i + (\varepsilon/3)B_X)_{i=1}^k $$
is a covering of ${\Bbb S}$. The thesis of the lemma and a trivial estimation of the Euclidean norm on $E$ gives that for some some $i \in \{ 1,\dots,k \}$ the set set $A_i + (\varepsilon/2)B_X$ contains ${\Bbb S} \cap F$ with $F$ an $n$-dimensional subspace of $E$. Therefore, $A_i + \varepsilon B_X$ contains $S_X \cap F$.\findemo\\

Now, we can give a satisfactory estimation of the radius of arbitrarily dimensional balls in the sets of a finite decomposition of $B_X$ into closed convex sets.

\begin{coro}\label{main_coro}
Let $n, k \in {\Bbb N}$ and $\lambda <1 $. There exists ${\mathcal N}^*(n, k, \lambda)$ such that for all the normed spaces $X$ of dimension $N \geq {\mathcal N}^*(n, k,\lambda)$ (that applies for $N=\infty$ in particular) and $(A_i)_{i=1}^k$ a covering of $B_X$ by closed convex bounded sets, there is some $i \in \{1, \dots, k\}$ such that $A_i$ contains an $n$-dimensional ball of radius $\lambda$.
\end{coro}

Note that the Lyusternik-Shnirel’man theorem says that {\it Borsuk's conjecture} is true for balls, being false in general \cite{book}. Corollary~\ref{main_coro} could be regarded somehow  as a strengthening of Lyusternik-Shnirel’man result.\\

In order to better exploit this result we introduce the following notion. The {\it asymptotic inradius} of a set $A \subset X$ is the number
$$ \rho(A):= \sup\{ \lambda>0: A~ \mbox{contains arbitrarily high-dimensional balls of radius} ~\lambda \}. $$

We have the following account of properties of $\rho$.

\begin{prop}\label{rules}
Let $X$ be an infinite-dimensional Banach space. Then
\begin{itemize}
\item[(a)] $\rho(B_X)=1$;
\item[(b)] if $A \subset B \subset X$ are convex closed bounded and $\lambda \geq 0$, then
$$ \rho(A) \leq \rho(B)  ~~\mbox{and}~~  \rho(\lambda A) = \lambda \, \rho(A) ;  $$
\item[(c)] if $A, A_1, \dots, A_k \subset X$ are convex closed and bounded such that $A \subset \bigcup_{j=1}^k A_j$, then
$$ \rho(A) \leq \max\{ \rho(A_j): 1\leq j \leq k\}.$$
\end{itemize}
\end{prop}

\noindent
{\bf Proof.} (a) and (b) are pretty evident. Note that (c) is Corollary \ref{main_coro} in case $A=B_X$. For the general case, take $r_A < \rho(A)$, $r < r_A$ and $\lambda=r/r_A$. Given $n \in {\Bbb N}$, find a ball $B$ of dimension  
${\mathcal N}^*(n, k,\lambda)$ inside $A$. Since $A_1, \dots, A_k$ is a cover of $B$, it reduces to Corollary \ref{main_coro} via homogeneity. One of the sets $A_1, \dots, A_k$ contains a ball of dimension $n$ and radius $r$.\findemo\\

Statement (c) of Proposition \ref{rules} applied to a finite convex decomposition of a convex set in an infinite-dimensional Banach space gives a ``principle of conservation'' the asymptotic inradius. We will state explicitly that as a corollary. 

\begin{coro}
Let $A, A_1, \dots, A_k \subset X$ be convex closed and bounded subsets of an infinite-dimensional Banach space such that $A = \bigcup_{j=1}^k A_j$, then there is some $j$ such that $\rho(A_j)=\rho(A)$.
\end{coro}

Despite the asymptotic inradius behave like a measure of noncompactness with respect to homothecies, translations and finite decompositions (into convex closed sets), it fails with respect to addition of sets.

\begin{ejem}
Consider the convex closed bounded set in $\ell_1$
$$ A= \{ (x_n) \in B_{\ell_1}: x_n \geq 0 \}.$$
Then $\rho(A)=\rho(-A)=0$, but $A+(-A) = B_{\ell_1}$.
\end{ejem}

\section{Covering $B_X$ by translates of a single set} 

In order to obtain a  result valid for convex non balanced sets we need to remove the symmetry in the well known Mazur's lemma \cite{LT}.

\begin{prop}
Let $X$ be a Banach space of infinite dimension, $A \subset X$ a bounded convex set that has $0$ as an interior point and let $\Lambda$ its Minkowski functional. If $F \subset X$ is a finite dimensional subspace and $\delta>0$ then there is a finite codimensional subspace $Y \subset X$ and a projection $P$ defined on $Z=F+Y$ onto $F$ with kernel $Y$ such that for every $z \in Z$ we have
$$ \Lambda (P(z)) \leq (1+\delta) \Lambda(z). $$
\end{prop}

\noindent
{\bf Proof.}  Let $\partial A$ be the boundary of $A$, that is the set $\{ \Lambda =1\}$. Note that $\Lambda$ is sublinear and Lipschitz with respect to the norm. The set $F \cap \partial A$ is compact so there are finitely many points $x_i \in F \cap \partial A$, with $1 \leq i \leq k$, such that for any $x \in F \cap \partial A$ there is $i$ such that $\Lambda(x_i-x) < \delta'$ where $\delta'=1-(1+\delta)^{-1}$.
By the Hahn-Banach theorem there are functionals $x_i^*$ such that $x_i^* \leq \Lambda$ and $x_i^*(x_i)=1$. Let $Y = \bigcap_{i=1}^k \ker x_i^*$. Assume that $x \in  F \cap \partial A$ and $y \in Y$. If $i$ is the suitable index then 
$$ \Lambda(x+y) \geq  \Lambda(x_i +y) - \Lambda(x_i-x)  \geq x_i^*(x_i+y) -\delta' = 1-\delta' = (1+\delta)^{-1} \Lambda(x) .$$
As the end members of the inequality are homogeneous then 
$$\Lambda(x) \leq (1+\delta) \Lambda(x+y)$$ 
whenever $x \in F$ and $y \in Y$. That also implies that the sum $F+Y$ is direct, so the projection $P$ is well defined.\findemo\\

We have the following general result.

\begin{theo}\label{teo_H}
Let $X$ be an infinite dimensional Banach space, $A \subset X$ a convex subset  and assume that $0 \in A$ and there is a finite subset $O \subset X$ such that $B_X \subset A+O$. Then for every $\varepsilon >0$ there is a finite codimensional subspace $Y \subset X$ and a finite subset $H \subset A$ such that 
$$ B_Y \subset (A - H) + \varepsilon B_X .$$
\end{theo}

\noindent
{\bf Proof.} 
Let $F$ be the finite-dimensional subspace generated by $O$, so we have
$B_X \subset A +F$. 
By Lemma \ref{tech_2}, $A$ has nonempty interior, so we may assume $0 \in A^{\mathrm{o}}$ without loss of generality.
Let $\Lambda$ the Minkowski functional of  $A$ and 
take $Y$ as in the proposition with $\delta>0$ such that
$$ (1+\delta) A \subset A + (\varepsilon/3)B_X $$
Consider the set $C =F \cap A  $
which is relatively compact, so we can find a finite subset $H \subset A$ such that 
$$C \subset H + (\varepsilon/3)B_X.$$
Now assume that $y \in B_Y$ there are $x \in D$ and $u \in F$ such that $y=x-u$. Clearly $x \in Z=Y+F$ so the projection $P$ is defined for $x$. We have 
$$ 0=P(y)= P(x) - P(u) = P(x) - u. $$
Observe that $u = P(x) \in (1+\delta)A$, and also $u \in F$. 
Since there is $v \in H$ such that $\|u-v\|<2\varepsilon/3$ and $z \in A$ with $\|x-z\| \leq \varepsilon/3$ we have
$$ y \in (A- H) + \varepsilon B_X.$$
as wanted.\findemo\\

The set $A$ satisfying the thesis of Theorem \ref{teo_H} may be smaller that the unit ball. Indeed, let $A \subset {\Bbb R}^2$ be a regular hexagon of radius ($=$ side) $1/\sqrt{3}$. If $H$ is the set of its vertices, then $B_{{\Bbb R}^2} \subset A-H$.  
 
\begin{coro}\label{coro_half}
Let $X$ be a  Banach space, $A \subset X$ a balanced bounded convex closed subset  and assume that there is a finite subset $O \subset X$ such that $B_X \subset A+O$. Then $A$ contains a finite-codimensional ball of radius $\lambda$ for every $\lambda<1/2$.
\end{coro}

\noindent
{\bf Proof.} 
Given $\varepsilon >0$, Theorem \ref{teo_H} provides a subset $H \subset A$ such that $B_X \subset (A-H)+\varepsilon B_X$. Note that for $A$ balanced we have  $A-H \subset A-A = 2A$. The final result is obtained by using Lemma~\ref{tech_1}.\findemo

\begin{coro}
Let $X$ be a  Banach space, $A \subset X$ a balanced bounded convex closed subset with nonempty interior. If there is a compact subset $K \subset X$ such that $B_X \subset A+K$, then $A$ contains a finite-codimensional ball of radius $\lambda$ for every $\lambda<1/2$.
\end{coro}

\noindent
{\bf Proof.} 
For every $\varepsilon >0$ there is finite subset $O \subset K$ such that $K \subset O + \varepsilon B_X$. Then apply Corollary \ref{coro_half} together Lemma \ref{tech_1}.\findemo\\

Note the similarity with Fredholm theory. However, here we are summing sets instead of operators.\\

We do not know how optimal is the constant $1/2$ obtained above. In the case of Hilbert space it can be improved assuming in addition that the covering set $A$ is contained in the unit ball.

\begin{prop}
Let $X$ be a Hilbert space, $A \subset B_X$ a  convex subset  and assume that there is $K \subset X$ compact such that $B_X \subset A+K$. Then, for every $\varepsilon >0$ there is a finite codimensional subspace $Y \subset X$ such that 
$$ B_Y \subset A + \varepsilon B_X .$$
\end{prop}

\noindent
{\bf Proof.} Given $\varepsilon>0$, by the uniform convexity of the norm, we may find $\delta>0$ such that if $x,y \in B_X$ and the segment $[x,y] $ is disjoint with $(1-\delta)B_X^{\mathrm{o}}$, then $\|x-y\| \leq \varepsilon$.
The compact $K$ is contained by the subspace generated by the first $n$ elements of a Hilbert basis except a perturbation of norm less than $\delta$. Let $Y$ be the subspace generated by the basis elements starting from $n+1$ on and let $P$ be the norm one projection onto $Y$. We have 
$B_Y \subset P(A) + \delta B_X$. 
Consider any $y \in S_Y$. Then, there is $z \in P(A)$ such that $\|z-y\| \leq \delta$ and $x \in A$ such that $P(x)=z$. For any $\lambda \in [0,1]$ we have
$$  P( \lambda y + (1 - \lambda) x )  =  \lambda y + (1 - \lambda) z .$$
Now, note that 
$$\| ( \lambda y + (1 - \lambda) z)-y \| = (1-\lambda)\|z-y\| \leq \delta, $$
and thus we  have $\|  \lambda y + (1 - \lambda) z\| \geq 1-\delta$.
That implies $\| \lambda y + (1 - \lambda) x \| \geq 1-\delta$ for every $\lambda \in [0,1]$.
The choice of $\delta$ implies that $\|x-y\| \leq \varepsilon$ and therefore $y \in A +  \varepsilon B_X$.
Since $y \in S_Y$ was arbitrary, we get $S_Y \subset A +  \varepsilon B_X$. Finally, by convexity we deduce $B_Y \subset A + \varepsilon B_X$.\findemo\\

The $\varepsilon$-perturbation cannot be removed as the following example in $\ell_2$ shows
$$ A = \{ (x_n) \in B_{\ell_2}:  \sum_{n=1}^\infty (1-2^{-n})^{-2} \,  x_n^2  \leq 1 \} $$
$$ \mbox{and}~~ K=  \{ (x_n) \in B_{\ell_2}: | x_n|  \leq 2^{-n} \} . $$

\begin{rema}
Continuous functions on the Hilbert sphere could not be almost constant on any infinite dimensional subspace after the negative solution to the distortion problem by Odell and Schlumprecht \cite{OS}, therefore a similar result to Corollary \ref{main_coro} with infinite dimensional subspaces seems to be not feasible, see \cite{Milman2} for a discussion in terms of the infinite spectrum.
\end{rema}

\section*{Acknowledgements}

My search for large finite codimensional balls actually started by thinking around the {\it Gorelik principle} and {\it asymptotic uniform smoothness} during my stay in Besançon in Spring 2018. I am indebted to Gilles Lancien for his hospitality and the long talks around the whiteboard. I wish to express my gratitude to Vitali Milman for providing me with valuable information on his results around Dvoretzky's theorem published in the 70's in Russian journals, mostly inaccessible to me. I am also grateful to Bernardo González Merino for stimulating discussions about estimates of the radius of Euclidean balls in low dimensions.

\vspace{1cm}

\begin{flushright}
Departamento de Matem\'aticas\\ Universidad de Murcia\\
Campus de Espinardo\\ 30100 Espinardo, Murcia, SPAIN\\
E-mail: matias@um.es
\end{flushright}

\end{document}